\documentclass[a4paper,11pt,oneside]{article}
\usepackage{amsmath,amssymb,amsthm,natbib,epsf}

\newcommand{\X}{\mathbf{X}}
\newcommand{\Y}{\mathbf {Y}}

\newcommand{\un}{\mathbf{1}}
\newcommand{\Nset}{\mathbb{N}}

\newcommand{\Rset}{\mathbb{R}}

\renewcommand{\epsilon}{\varepsilon}

\def\prob{\mathbb{P}}
\def\probs{\overline{\mathbb{P}}}
\def\esp{\mathbb{E}}
\def\esps{\overline{\mathbb{E}}}

\def\denss{\overline{p}}
\def\un{\mathbf{1}}

\def\Xb{\mathbf{X}}

\def\Yb{\mathbf{Y}}

\newcommand{\gtnu}[2]{\mbox{$g_{\theta}(#1|#2)$}}

\newcounter{hyp}
\newenvironment{hyp}[1]{\refstepcounter{hyp}
  \begin{list}{}{\leftmargin=11mm\labelwidth=\leftmargin}
  \item[({\bf A\arabic{hyp}})] \label{#1}}{\end{list}}

\newcounter{hypI}

\newcommand{\refhyp}[1]{{\rm (}{\bf A\ref{#1}}{\rm )}}

\newcommand{\refhyps}[2]
   {{\rm (}{\bf A\ref{#1}}{\rm )}--{\rm (}{\bf A\ref{#2}}{\rm )}}

\newtheorem{thm}{Theorem}
\newtheorem{lemma}{Lemma}

\newtheorem{prop}{Proposition}

\theoremstyle{remark}

\newcommand{\nat}{\nabla_\theta}

\renewenvironment{abstract}%
 {\bigskip \noindent \small \bf {ABSTRACT.}}%
 {}

\begin{document}

\def\thefootnote{\fnsymbol{footnote}}

\thispagestyle{empty}

\bigskip

\noindent
  {\bf\Large  Non singularity of the asymptotic Fisher information matrix in hidden Markov models}\\
  \ \\
  RANDAL DOUC\footnote{Randal DOUC, CMAP, \'Ecole Polytechnique, Route de Saclay, 9128 Palaiseau Cedex. FRANCE. {\sf douc@cmapx.polytechnique.fr}}\\
{\em \'Ecole Polytechnique}
 
\bigskip

\begin{abstract}
  In this paper, we consider a parametric hidden Markov model where the hidden state space is non necessarily finite. We provide a necessary and sufficient condition for the invertibility of the limiting Fisher information matrix. 

\end{abstract}
 

\bigskip
\noindent {\it AMS classification codes:}
  {\small Primary 62M09. Secondary 62F12.}

\newpage

\section{Introduction}
Hidden Markov Models (HMMs) form a wide class of discrete-time stochastic
processes, used in different areas such as speech recognition (\cite{juang:rabiner:1991}), neurophysiology (\cite{fredkin:rice:1987}), biology (\cite{churchill:1989}), and time series analysis (\cite{dejong:shephard:1995}, \cite{chan:ledolter:1995}; see \cite{macdonald:zucchini:1997} and the references therein).

The motivation for finding conditions implying the non singularity of the limiting Fisher information matrix is intrinsequely linked with the asymptotic properties of the maximum likelihood estimator (MLE) in those models. Since the way of dealing with consistency and asympotic normality of the MLE is not yet unified for hidden Markov models, we first give a quick overview of the different types of techniques that were developped aroud the MLE properties.  

Most works on maximum likelihood estimation in such models have focused on
iterative numerical methods, suitable for approximating the maximum 
likelihood estimator. By contrast, the statistical issues regarding the asymptotic properties of the
maximum likelihood estimator itself have been largely ignored
until recently. \cite{baum:petrie:1966} have shown the consistency and
asymptotic normality of the maximum likelihood estimator in the particular case
where both the observed and the latent variables take only finitely many
values. These results have been extended recently in series of papers by \cite{leroux:1992}, \cite{bickel:ritov:ryden:1998}, \cite{jensen:petersen:1999} and \cite{douc:moulines:ryden:2004}. The latter authors generalize the method followed by \cite{jensen:petersen:1999} to swiching autoregressive models associated to a possibly non finite hidden state space, the observations belonging to a general topological space. Their method put the consistency and the asymptotic normality in a common framework where some ``stationary approximation'' is performed under uniform ergodicity of the hidden Markov chain. This stringent assumption seems hard to check in a non compact state space. Nevertheless, up to our best knowledge,  the assumptions used in \cite{douc:moulines:ryden:2004} are the weakest known in hidden Markov models literature for proving these asymptotic results, even if the extension to a non compact state space is still an open question.    

Another approach was initiated by \cite{legland:mevel:2000}. They independently developed
a different technique to prove the consistency and the asymptotic normality of the MLE (Mevel 1997) for hidden Markov models with finite
hidden state space. The work of \cite{legland:mevel:2000} later generalised to a non finite state space by \cite{douc:matias:2002} is based on the remark
that the likelihood can be expressed as an additive function of an ``extended
Markov chain''. They show that under appropriate conditions, this extended chain is in some sense geometrically
ergodic once again under the assumption of uniform ergodicity for the hidden Markov chain. Nevertheless, even if the ergodicity of the extended Markov chain may be of independent interest, the assumptions used in this approach are stronger than in \cite{douc:moulines:ryden:2004}. 

However, in all these papers, the asymptotic normality of the MLE is derived from the consistency property thanks to the non singularity of the limiting Fisher information matrix. Indeed, whatever approach is considered (``extended Markov chain'' method or ``stationary approximation'' approach), the asymptotic normality of the MLE is obtained through a Taylor expansion of the gradient of the loglikelihood around the true value of the parameter. The given equation is then transformed by inversion of the Fisher information matrix associated to $N$ observations so as to isolate the quantity of interest $\sqrt{N}(\theta_{MV}-\theta^*)$ (where $\theta_{MV}$ is the maximum likelihood estimator and $\theta^*$ the true value). The last step consists in considering the asymptotic behavior of the obtained equation as the number of observations grows to infinity and in particular, a crucial feature is that the normalised information matrix should converge to a non singular matrix. 

Unfortunately, the asymptotic Fisher information matrix $I(\theta)=-\esp_\theta\left( \nat^2 \log p_\theta(Y_1|\Y_{-\infty}^0)\right)$ is the expectation of a quantity which does depend on all the previous observations. Under this form, the non singularity of this matrix is hardly readable. The aim of this paper is to show that this non singularity is equivalent to the non singularity of some $I_{Y_{1:n}}(\theta)=-\esp_\theta\left( \nat^2 \log p_\theta(Y_{1:n})\right)$ where here, the expectation only concerns a finite number of observations. One thus might expect that non singularity of $I_{Y_{1:n}}(\theta)$ is easier to check than the one of $I(\theta)$. This is a simple result but we expect that it helps for checking such intractable non singularity assumption. The rest of the paper is organised as follows: we introduce the model and the assumptions in Section 2. In Section 3, after recalling some properties of the MLE, we state and prove the main result of the paper using a technical proposition. Finally, Section 4 is devoted to the proof of this technical proposition.

\section{Model and assumptions}               
In the following, the assumptions on the model and the description of the asymptotic results concerning the MLE directly derive from the paper of \cite{douc:moulines:ryden:2004}. Let $\{X_n\}_{n=0}^\infty$ be a Markov Chain on $(\X,{\mathcal B}(\X))$. We denote by $\{Q_\theta(x,\mathcal {A}), x \in \X, \mathcal{A} 
\in 
\mathcal{B(}\X)\}$, the
Markov transition kernel of the chain. We also 
let 
$\{Y_n\}_{n=0}^\infty$ be a sequence of
random variables in $(\Y,{\mathcal B}(\Y))$, 
such
that, conditional on $\{X_n\}_{n=0}^\infty$, $\{Y_n\}_{n=0}^\infty$ is a
sequence of conditionally independent random variables, $Y_n$ with 
conditional
density $\gtnu{y}{X_n}$ with respect to some $\sigma$-finite measure $\nu$ 
on
the Borel $\sigma$-field $\mathcal {B(\Y)}$. Usually, $\X$ and ${\Y}$ 
are 
subsets of $\Rset^{s}$ and
$\Rset^{t}$ respectively, but they may also be higher dimensional spaces.
Moreover, both $Q_\theta$ and $g_\theta$ depend on a parameter $\theta$ in
$\Theta$, where $\Theta$ is a compact subset of $\Rset^{p}$. The true parameter
value will be denoted $\theta^\ast$ and is assumed to be in the interior of 
$\Theta$.

Assume that for any $x \in \X$, $Q_{\theta}(x, \cdot)$ has a density $q_{\theta}(x,\cdot)$ with
respect to the same $\sigma$-finite dominating measure $\mu$ on $\mathcal
{X}$. For any $k \geq 1$, the density of $Q_{\theta}^k(x, \cdot)$ with respect to $\mu$ is denoted by $q_{\theta}^k(x, \cdot)$. 
In the following, for $m \geq n$, denote $\Yb_n^m$ the family
of random variables $(Y_n, \ldots, Y_m)$.  Moreover, for any measurable function $f$ on
$\mathcal{(\X,B}(\X),\mu)$, denote ess sup($f)=\inf\{M\geq 0, \mu(\{M<
|f|\})=0\}$ and if $f$ is non-negative, ess inf($f)=\sup\{M\geq 0,
\mu(\{M > f\})=0\}$ (with obvious conventions if those sets are empty).
By convention, we simply write $\sup$ (resp. $\inf$) instead of $\mbox{ess sup}$ (resp. $\mbox{ess inf}$).  

We denote by $\pi_\theta$ the stationary distribution of the kernel $Q_\theta$ when it exists. Let $(Z_n=(X_n,Y_n))_{n \geq 0}$ be the Markov chain on $((\X \times \Y)^{\Nset}, ({\mathcal B}(\X)\otimes {\mathcal B}(\Y))^{\otimes \Nset})$ of transition kernel $P_\theta$ defined by 
$$
P_{\theta}(z,A):= \iint \un_{A}(x',y') Q_\theta(x,dx') g_{\theta}(y'|x') \nu(dy').
$$
$\prob_{\theta,\lambda}$ (resp. $\esp_{\theta,\lambda}$) denotes the probability (resp. expectation) induced by the Markov chain $(Z_n)_{n \geq 0}$ of transition kernel $P_\theta$ and initial distribution $\lambda(dx)g_\theta(y|x)\nu(dy)$. By convention, we simply write $\prob_{\theta,x}:=\prob_{\theta,\delta_{x}}$ and $\probs_{\theta}:=\prob_{\theta,\pi_\theta}$ and $\esp_{\theta,x}$ and $\esps_{\theta}$ will be the associated expectation. Moreover, $\pi_\theta$ will also denote the density of the stationary distribution with respect to $\mu$. 

In this paper, $\|\cdot\|$ denotes the $L_2$ norm in $\Rset^p$, i.e. for any $\varphi \in \Rset^p$, $\|\varphi\|=\sqrt{\varphi^T\varphi}$. By abuse of notation, $\|\cdot\|$ will also denote the associated $L^2$-norm in the space of symmetric matrices in $\Rset^p \times \Rset^p$, i.e. for any real $p\times p$ matrix $J$, $\|J\|:=\sup_{\varphi, \|\varphi\|=1}|\varphi^T J \varphi|$. For all bounded measurable function $f$, define $\|f\|_\infty:=\sup_{x \in \X}|f(x)|$. $C$ will denote unspecified finite constant which may take different
values upon each appearance. In the following, we will use $p_\theta$ as a generic symbol for density. When this density explicitely depends on $\pi_\theta$, we stress it by writing $\bar p_\theta$ instead of $p_\theta$. In case of ambiguity, we will define precisely the density $p_\theta$.     

\noindent \textbf{Consistency assumptions.} I recall the assumptions, used in \cite{douc:moulines:ryden:2004} to obtain consistency for a switching autoregressive model. Since we consider here a hidden Markov model, we adapt the statement of their assumptions. 

\begin{hyp}{hyp:uniform-ergodicity}
\begin{itemize}
\item[(a)] $0<\sigma_-:=\inf_{\theta\in\Theta}\inf_{x,x'\in\X}
  q_\theta (x,x')$
and $\sigma_+ :=
 \sup_{\theta\in\Theta}\sup_{x,x'\in\X} q_\theta(x,x')<\infty$.
\item[(b)] For all $y\in\Y$,
$0<\inf_{\theta\in\Theta}\int_\X g_\theta(y|x)\,\mu(dx)$
and
$\sup_{\theta\in\Theta}\int_\X g_\theta(y|x)$ $\mu(dx)<\infty$.
\end{itemize}
\end{hyp}
\begin{hyp}{hyp:moment-condition}
$b_+:= \sup_\theta \sup_{y,x} g_\theta(y| x)
  <\infty$
and $\esps_{\theta^\ast}(|\!\log b_-(Y_1)|)<\infty$,
where $b_-(y):= \inf_\theta
  \int_{\X}g_\theta(y|x)\,\mu(dx)$.
\end{hyp}

\begin{hyp}{hyp:identifiability}
$\theta=\theta^\ast$ if and only if $\probs_{\theta}^Y=\probs_{\theta^\ast}^Y$, where $\probs_{\theta}^Y$ is the trace of $\probs_{\theta}$ on
$\{\Y^\Nset,\mathcal{B}(\Y)^{\otimes \Nset}\}$, that is the distribution
of $\{Y_k\}$.
\end{hyp}

\noindent \textbf{Asymptotic normality assumptions.} Some additional assumptions are needed for the asymptotic normality of the MLE. We will assume that there exists a positive real $\delta$ such
that on
$G:=\{\theta\in\Theta: \|\theta-\theta^\ast\|<\delta\}$,
the following conditions hold.

\begin{hyp}{hyp:differentiability}
For all $x,x'\in\X$ and $y'\in \Y$, the functions
  $\theta \mapsto q_\theta(x,x')$ and $\theta\mapsto g_\theta(y'|x')$
  are twice continuously differentiable on $G$.
\end{hyp}

\begin{hyp}{hyp:Derivative-Hessian}
\begin{itemize}
\item[(a)]
 $\sup_{\theta\in G}\sup_{x,x'}\|\nat\log q_\theta(x,x')\|<\infty$,
 $\sup_{\theta\in G}\sup_{x,x'}\|\nat^2\log q_\theta(x,x')\|<\infty$.
\item[(b)] $\esps_{\theta^\ast}[\sup_{\theta\in G}\sup_x
  \|\nat\log g_\theta(Y_1|x)\|^2]<\infty$,\!
 $\esps_{\theta^\ast}[\sup_{\theta\in G}\sup_x
 \|\nat^2 \log g_\theta(Y_1|x)\|]$ $<\infty$.
\end{itemize}
\end{hyp}

\begin{hyp}{hyp:UniformBound}
\begin{itemize}
 \item [(a)] For $\nu$-almost all $y$ in
  $\Y$ there exists a
  function $f_y:\,\X\to\Rset^{+}$ in $L^1(\mu)$ such
  that $\sup_{\theta\in G} g_\theta(y|x)\leq f_{y}(x)$.
 \item[(b)] For $\mu $-almost all $x\in\X $, there exist functions
   $f^1_{x}:\,\Y\to\Rset^{+}$ and
   $f^2_{x}:\,\Y\to\Rset^{+}$ in $L^1(\nu)$
   such that
   $\|\nat g_\theta(y|x)\|\leq f^1_{x}(y)$
   and $\|\nat^2 g_\theta(y|x)\|\leq f^2_{x}(y)$
   for all $\theta\in G$.
\end{itemize}
\end{hyp}

\section{Main result}
We first recall the results of consistency and asymptotic normality of the MLE obtained by \cite{douc:moulines:ryden:2004}. Write $\hat \theta_{n,x_0}$ the maximum likelihood estimator associated to the initial condition $X_0 \sim \delta_{x_0}$ where $\delta_{x_0}$ is the dirac mass centered in $x_0$. 
\begin{thm}\label{theo:Consistency-MLE-2}
Assume \refhyps{hyp:uniform-ergodicity}{hyp:identifiability}. Then, for
any $x_0\in\X$, 
$$
\lim_{n\to\infty}\hat \theta_{n,x_0}=\theta^\ast \quad \probs_{\theta^\ast}-a.s.
$$
\end{thm}
Define the stationary density $\bar p_\theta$ and the Fisher information matrix $I_{\Yb_1^n}(\theta)$ associated to $n$ observation by:
\begin{align*}
&\denss_\theta(\Yb_1^n) :=\idotsint \pi_\theta(dx_1)\prod_{i=1}^{n-1} q_\theta(x_i,x_{i+1})\prod_{i=1}^n g_\theta(Y_i|x_i) \mu^{\otimes n}(dx_1^n),\\
&I_{\Yb_1^n}(\theta) :=-\esps_{\theta}(\nat^2 \log \denss_{\theta}(\Yb_{1}^n)). 
\end{align*}  
\cite{douc:moulines:ryden:2004} (Section 6.2) proved that $\lim_{n \to \infty} -\esps_{\theta}(\nat^2 \log \denss_{\theta}(\Yb_{1}^n))/n$ exists and denoting by  $I(\theta)$ this limit, they obtained that the asymptotic Fisher information matrix $I(\theta)$ may also be written as 
\begin{equation}
  \label{eq:expressionFisher}
 I(\theta)=\lim_{n \to \infty} -\esps_{\theta}(\nat^2 \log \denss_{\theta}(\Yb_{1}^n))/n= \lim_{m \to \infty} \esps_{\theta}(-\nat^2 \log \denss_{\theta}(Y_1|\Yb_{-m}^0)).
\end{equation}
We recall the asymptotic normality of the MLE obtained by \cite{douc:moulines:ryden:2004}. 
\begin{thm}\label{thm:tcl}
Assume \refhyps{hyp:uniform-ergodicity}{hyp:UniformBound}. Then, provided that $I(\theta^\ast)$ is non singular, we have for any $x_0\in\X$,
$$
\sqrt{n}(\hat \theta_{n,x_0}-\theta^*)\to\mathcal{N}(0,I(\theta^\ast)^{-1})
 \quad \mbox{$\probs_{\theta^\ast}$-\textrm{weakly}}.
$$
\end{thm}

We may now state the main result of the paper which links the asymptotic Fisher information matrix with the stationary information matrix associated to a finite number of observations. 
\begin{thm}
\label{thm:mainresult}
  Assume \refhyps{hyp:uniform-ergodicity}{hyp:UniformBound}. Then, $I(\theta^*)$ is non singular if and only if there exists $ n \geq 1$, such that $I_{\Yb_1^n}(\theta^*)$ is non singular.
\end{thm}
Note that this theorem holds under the same assumptions as in Theorem \ref{thm:tcl}. Of course, since the true parameter is not known, the non singularity of $I_{\Yb_1^n}(\theta)$ should be checked for all $\theta$. More precisely, if the stationary distribution $\pi_\theta$ of the hidden Markov chain is sufficiently known so that for any parameter $\theta$, the Fisher information matrix $I_{\Yb_1^n}(\theta)$ is shown to be non singular (for some $n$), then the sufficient condition of Theorem \ref{thm:mainresult} ensures that the asymptotic Fisher information matrix is non singular. Thus, the MLE is asymptotically normal by applying Theorem \ref{thm:tcl}. Before proving the necessary and sufficient condition of Theorem \ref{thm:mainresult}, we need a technical proposition about some asymptotics of the Fisher information matrix. Let 
\begin{align*}
  &I_{\Yb_1^n| \Yb_{-k-m}^{-k}}({\theta}):=-\esps_{{\theta^*}}(\nat^2 \log \denss_{{\theta^*}}(\Yb_1^n| \Yb_{-k-m}^{-k})).
\end{align*}  
\begin{prop} \label{prop:stationFisher}
  Assume \refhyps{hyp:uniform-ergodicity}{hyp:UniformBound}. Then, for all $n \geq 1$, 
$$
\lim_{k\to \infty } \sup_{m}\|I_{\Yb_{1}^n| \Yb_{-k-m}^{-k}}({\theta^*})-I_{\Yb_{1}^n}({\theta^*})\|=0.
$$ 
\end{prop}
While this result seems intuitive since the ergodicity of $(X_k)$ implies the asymptotic independence of $(\Yb_1^n)$ (where $n$ is fixed) wrt $\sigma\{ \Yb_{-l}^{-k}; l\geq k\}$, the rigourous proof of this proposition is rather technical and is thus postponed to Section 4. Using this proposition, Theorem \ref{thm:mainresult} may now be proved with elementary arguments. 

\begin{proof} (Theorem \ref{thm:mainresult})
By Eq. (\ref{eq:expressionFisher}), 
$$
\mathrm{det} I(\theta^*)=\mathrm{det}\left(\lim_n \frac{I_{\Yb_1^n}(\theta^*)}{n}\right)=\lim_n \mathrm{det}\left(\frac{I_{\Yb_1^n}(\theta^*)}{n}\right). 
$$ 
And thus if for all $n$, $I_{\Yb_1^n}(\theta^*)$ is singular then $I(\theta^*)$ is singular. Now, assume that $I(\theta^*)$ is singular. Fix some $n\geq 1$ and let $k \geq n$. 
By stationarity of the sequence $(Y_i)_{i \in {\mathbb Z}}$ under $\probs_{\theta}$ and elementary properties of the Fisher information matrix, 

\begin{align*}
&kI({\theta^*})=k \lim_{m \to \infty} \esps_{{\theta^*}}(-\nat^2 \log \denss_{{\theta^*}}(Y_1|\Yb_{-m}^0)),\\
&\quad = \lim_{m \to \infty} \sum_{i=1}^k  \esps_{{\theta^*}}(-\nat^2 \log \denss_{{\theta^*}}(Y_{i}|\Yb_{-m+i-1}^{i-1})),\\
&\quad =\lim_{m \to \infty}  \esps_{{\theta^*}}(-\nat^2 \log \denss_{{\theta^*}}(\Yb_1^{k}|\Yb_{-m}^0)),\\
&\quad =\lim_{m \to \infty}  \left[\esps_{{\theta^*}}(-\nat^2 \log \denss_{{\theta^*}}(\Yb_{k-n+1}^{k}|\Yb_{-m}^0 ))+ \esps_{{\theta^*}}(-\nat^2 \log \denss_{{\theta^*}}(\Yb_1^{k-n}|\Yb_{k-n+1}^{k}, \ \Yb_{-m}^0 ))\right],\\
 &\quad \geq\limsup_{m \to \infty} \esps_{{\theta^*}}(-\nat^2 \log \denss_{{\theta^*}}(\Yb_{k-n+1}^{k}|\Yb_{-m}^0))= \limsup_{m \to \infty} \esps_{{\theta^*}}(-\nat^2 \log \denss_{{\theta^*}}(\Yb_1^{n}| \Yb_{-k+n-m}^{-k+n})).
\end{align*}

\noindent Let $\varphi \in \Rset^p$ such that $I({\theta^*})\varphi=0$. Then, by the above inequality, for all $k\geq n$,   
$$
\varphi^T \limsup_{m \to \infty} \esps_{{\theta^*}}(-\nat^2 \log \denss_{{\theta^*}}(\Yb_1^n| \Yb_{-k+n-m}^{-k+n}))\varphi=0. 
$$
Now, letting $k\rightarrow \infty$, and using Proposition~\ref{prop:stationFisher}, we get $I_{\Yb_1^n}({\theta^*})\varphi =0$.  Thus, $I_{\Yb_1^n}$ is singular for any $n\geq 1$. The proof is completed. 
  
\end{proof}

\section{Proof of Proposition \ref{prop:stationFisher}}

\noindent \textbf{Regularity of the stationary distribution.} We first check that $\theta \mapsto \pi_\theta(f)$ is twice differentiable et obtain a closed form expression of $\nat \pi_\theta(f)$. By \refhyp{hyp:uniform-ergodicity}, the Markov chain $\{X_k\}_{k \geq 0}$ is uniformly ergodic and thus the Poisson equation associated to:
$$
V (x)-Q_{\theta^*}(x,V)=f(x)-\pi_{\theta^*}(f) \quad \mbox{ where } \quad \pi_{\theta^*}(|f|)<\infty,
$$ 
admits a unique solution that we denote by  $V_{\theta^*}$. Classically, since for all $(x,A)$ in $(\X,{\mathcal B}(\X))$, $Q_{\theta^*}(x,A)>\sigma_- \mu(A)$, we have
for all bounded measurable $f$, 
$$
|Q^k_\theta(x,f)-\pi_\theta(f)| \leq 2 \|f\|_\infty (1-\sigma_-)^k,
$$ 
and then
\begin{equation} \label{eq:majoration}
\forall x \in \X, \quad |V_{\theta^*}(x)|= \left|\sum_{k=0}^\infty \left(Q_{\theta^*}^k(x,f)-\pi_{\theta^*}(f)\right)\right| \leq 2\frac{\|f\|_\infty}{\sigma_-}. 
\end{equation}
Note that for all $\theta \in \Theta$ and for all bounded measurable function $f$,  
\begin{eqnarray}
\pi_\theta(f)-\pi_{\theta^*}(f)&=&\int \pi_\theta(dx)(f(x)-\pi_{\theta^*}(f))=\int \pi_\theta(dx)(V_{\theta^*}(x)-Q_{\theta^*}(x,V_{\theta^*}), \nonumber \\
&=&\pi_\theta \left( Q_\theta -Q_{ \theta^*}\right)V_{\theta^*}, \label{eq:differencepi}
\end{eqnarray}
where we have used that $\pi_{\theta}$ is the invariant probability measure for the transition kernel $Q_{\theta}$. Eq.~(\ref{eq:differencepi}) implies under\refhyps{hyp:uniform-ergodicity}{hyp:UniformBound} that 
 \begin{equation*} 
    \nat \pi_{\theta^*}(f)= \iint \pi_{{\theta^*}}(dx) \nat q_{\theta^*}(x,x')V_{\theta^*}(x')\mu(dx')= \esps_{{\theta^*}}\left\{\nat \log q_{\theta^*}(X_0,X_1)V_{\theta^*}(X_1)\right\}.\\
  \end{equation*}
Combining the uniform ergodicity ($\sigma_-=\inf_{\theta} \inf_{x,x'} q_{\theta}(x,x')>0$) with \refhyp{hyp:differentiability}, the conditions 1-3 of \cite{heidergott:hordijk:2003} are trivially satisfied and Theorem 4 of \cite{heidergott:hordijk:2003} ensures that ${\theta} \mapsto \pi_{\theta}(f)$ is twice differentiable. Moreover, applying Eq. (\ref{eq:closedfderiv}) to the bounded function $f(\cdot)=q_{\theta^*}(\cdot,x)$ where $x$ is a fixed point in $\X$ yields that ${\theta} \mapsto \pi_{\theta}(x)$ is twice differentiable at $\theta=\theta^*$ and
\begin{equation}\label{eq:closedfderiv}
\nat \pi_{\theta^*}(x)= \esps_{{\theta^*}}\left\{\nat \log q_{\theta^*}(X_0,X_1)\sum_{k=0}^\infty \left(q_{\theta^*}^{k+1}(X_1,x)-\pi_{\theta^*}(x)\right)\right\}.
\end{equation}

\noindent \textbf{Technical bounds.} We will now state and prove some technical bounds that will be useful for the proof of Proposition \ref{prop:stationFisher}. Lemma 9 in \cite{douc:moulines:ryden:2004} ensures the uniform forgetting of the initial distribution for the reverse a posteriori chain. It implies that for all $m\leq n,\ 0 \leq k \leq n-m$,
\begin{equation} \label{eq:oubliReverse}
\forall A \in {\mathcal B}(\X), \quad | \probs_\theta(X_{n-k} \in A|\Yb_{m}^n,X_n=x)-\probs_\theta(X_{n-k} \in A|\Yb_{m}^n, X_n=x')| \leq \rho^k,
\end{equation}
where $\rho:=1-{\sigma_-}/{\sigma_+}$. 

\begin{lemma} \label{lem:distribInit}
Assume \refhyps{hyp:uniform-ergodicity}{hyp:UniformBound}, then we have the following inequalities. 
  \begin{itemize}
\item[(i)]$\|\int \nat \pi_{\theta^*}(x)f(x)\mu(dx)\| \leq C \|f\|_\infty$ with $C:=2\frac{\sup_{\theta^*} \| \nat \log q_{\theta^*}(\cdot,\cdot)\|_\infty}{\sigma_-}$. 
\item[(ii)] For all $k\geq 0$, there exists a random variable $D(\Yb_{-\infty}^{-k})$ satisfying $\esps_{{\theta^*}}(D(\Yb_{-\infty}^{-k})^2) < \infty$ such that for all $m\geq 0$,   
$$
\left\|\int \nat \denss_{\theta^*}(x_{-k+1}|\Yb_{-k-m}^{-k}) f(x_{-k}) \mu(dx_{-k+1})\right\| \leq D(\Yb_{-\infty}^{-k}) \|f\|_\infty.
$$
Moreover,  for all $ k \geq 0$, $\esps_{{\theta^*}}(D(\Yb_{-\infty}^{-k})^2)= \esps_{{\theta^*}}(D(\Yb_{-\infty}^{0})^2)$. 
\end{itemize}
\end{lemma}
\begin{proof}
  Combining Eq (\ref{eq:majoration}) and Eq (\ref{eq:closedfderiv}) yields the first inequality. We will prove the second inequality with $k=0$. It is actually sufficient to bound $\sup_{x_1}|\nat \log \denss_{\theta^*}(x_1|\Yb_{-m}^0)|$. This can be done by using the Fisher identity: 
  \begin{align*}
&    \nat \log \denss_{\theta^*}(X_1|\Yb_{-m}^0)= \esps_{\theta^*} (\nat \log\denss_{\theta^*} (\Xb_{-m}^1,\Yb_{-m}^0)|X_1,\Yb_{-m}^0)-\esps_{\theta^*} (\nat \log\denss_{\theta^*} (\Xb_{-m}^1,\Yb_{-m}^0)|\Yb_{-m}^0),\\
&= \int \esps_{{\theta^*}} \left\{\left. \nat \log \pi_{\theta^*} (X_{-m})+\sum_{i=-m}^{-1} [\nat \log q_{\theta^*}(X_i,X_{i+1})+\nat \log g_{\theta^*}(Y_i|X_i)]\right|X_0,\Yb_{-m}^0\right\}\\
&\quad\quad [\probs_{\theta^*}(dX_{0}|X_1,\Yb_{-m}^0)-\probs_{\theta^*}(dX_{0}|\Yb_{-m}^0)]+2  [\|\nat \log q_{\theta^*}(\cdot,\cdot)\|_\infty +\|\nat \log g_{\theta^*}(Y_0|\cdot)\|_\infty ].\\
  \end{align*}
Thus, using Eq. (\ref{eq:oubliReverse}), we get
\begin{align*}
&  \sup_{x_1}\|\nat \log \denss_{\theta^*}(x_1|\Yb_{-m}^0)\| \\
&\quad \leq \|\nat \log \pi_{\theta^*}\|_\infty \rho^m+\sum_{i=-\infty}^0 \rho^{|i|}(1+\mathbf{1}_{i=0}) [\rho^{-1}\|\nat \log q_{\theta^*}(\cdot,\cdot)\|_\infty +\|\nat \log g_{\theta^*}(Y_i|\cdot)\|_\infty ],\\
 &\quad \leq  \|\nat \log \pi_{\theta^*}\|_\infty + \sum_{i=-\infty}^0 (1+\mathbf{1}_{i=0}) \rho^{|i|} [\rho^{-1}\|\nat \log q_{\theta^*}(\cdot,\cdot)\|_\infty +\|\nat \log g_{\theta^*}(Y_i|\cdot)\|_\infty ]= D(\Yb_{-\infty}^0).
\end{align*}
Using \refhyp{hyp:Derivative-Hessian}, it is straightforward that $\esps_{{\theta^*}}(D(\Yb_{-\infty}^{0})^2) < \infty$. Morever, by stationarity of $(Y_i)$ under $\probs_{\theta^*}$,  $\esps_{{\theta^*}}(D(\Yb_{-\infty}^{-k})^2)= \esps_{{\theta^*}}(D(\Yb_{-\infty}^{0})^2) $. The proof is completed. 
\end{proof}
Define $p_{\theta}(\Yb_1^n|X_{-k}=x):=\idotsint q^{k}_\theta(x,x_0) \prod_{i=1}^n q_{\theta}(x_{i-1},x_i) g_{\theta}(Y_i|x_i)\mu^{\otimes (n+1)}(dx_0^n)$.
\begin{lemma} \label{lem:techniq}
  Assume \refhyps{hyp:uniform-ergodicity}{hyp:UniformBound} and fix some $n \geq 1$. Then,  
  \begin{itemize}
  \item [(i)] For all $k \geq 0$, 
$$
\frac{\sup_{x,x'}\left| p_{\theta^*}(\Yb_1^n|X_{-k}=x)- p_{\theta^*}(\Yb_1^n|X_{-k}=x')\right|} {\inf_x p_{\theta^*}(\Yb_1^n| X_0=x)}\leq 2 (1-\sigma_-)^{k} \left(\frac{\sigma+}{\sigma_-}\right).
$$
\item  [(ii)] There exists $1>\lambda>1-\sigma_-$ and a random variable $E(\Yb_1^n)$ satisfying $\esp( E(\Yb_1^n)^2)<\infty$ such that for all $k \geq 0$, 
$$
\frac{\sup_{x,x'}\left\| \nat p_{\theta^*}(\Yb_1^n|X_{-k}=x)-\nat p_{\theta^*}(\Yb_1^n|X_{-k}=x')\right\|}{\inf_x p_{\theta^*}(\Yb_1^n| X_0=x)} \leq \lambda^k E(\Yb_1^n) . 
$$   
  \end{itemize}
\end{lemma}

\begin{proof}
Define $f_{\theta}(x):=p_{\theta}(\Yb_1^n|X_0=x)$. 
First, note that 
\begin{equation}
  \label{eq:majoVrais}
  \frac{\sup_x f_{\theta^*}(x)}{\inf_x f_{\theta^*}(x)} \leq  \frac{\sup_x \int p_{\theta^*}(\Yb_1^n|X_1=x_1)Q_{\theta^*}(x,dx_1)}{\inf_x \int p_{\theta^*}(\Yb_1^n|X_1=x_1)Q_{\theta^*}(x,dx_1)}\leq \frac{\sigma_+}{\sigma_-}. 
\end{equation}
 For all $x,x'$ in $\X$, 
  \begin{align*}
&    | p_{\theta^*}(\Yb_1^n|X_{-k}=x)-p_{\theta^*}(\Yb_1^n|X_{-k}=x')| =\left|\int f_{\theta^*}(x)(Q_{\theta^*}^k(x,dx_0)-Q_{\theta^*}^k(x',dx_0))\right|\\
& \quad \leq  2 (1-\sigma_-)^{k} \sup_x f_{\theta^*}(x),
  \end{align*}
where we have used the uniform ergodicity of the chain $(X_n)$. Combining with (\ref{eq:majoVrais}) completes the proof of (i). Now, note that  $p_{\theta}(\Yb_1^n|X_{-k}=x)=Q_{\theta}^{k}(x,f_{\theta})$ and write 
\begin{eqnarray*}
\nat p_{\theta^*}(\Yb_1^n|X_{-k}=x)&=&\nat \left. \left[Q_{\theta}^{k}(x,f_{\theta})\right]\right|_{\theta=\theta^*},\\
&=&\esp_{\theta^*,x}(\nat f_{\theta^*}(X_k))+  \esp_{\theta^*,x}\left\{f_{\theta^*}(X_k) \left( \nat \log  \prod_{i=1}^k q_{\theta^*}(X_{i-1},X_{i})\right)\right\},\\
&=&\esp_{\theta^*,x}(\nat f_{\theta^*}(X_k))+  \esp_{\theta^*,x}\left\{(f_{\theta^*}(X_k) -\pi_{\theta^*}(f_{\theta^*}))\left( \nat \log  \prod_{i=1}^k q_{\theta^*}(X_{i-1},X_{i})\right)\right\},\\
&=&\esp_{\theta^*,x}(\nat f_{\theta^*}(X_k))+ \sum_{i=1}^k \esp_{\theta^*,x}\left\{ \left(Q^{k-i}(X_i,f_{\theta^*})-\pi_{\theta^*}(f_{\theta^*})\right) \nat \log q_{\theta^*}(X_{i-1},X_{i})\right\}.
\end{eqnarray*}
Using that the chain $(X_i)$ is uniformly ergodic yields  
\begin{align*}
  &\left\| \nat p_{\theta^*}(\Yb_1^n|X_{-k}=x)- \nat p_{\theta^*}(\Yb_1^n|X_{-k}=x')\right\|\\
&\leq  2 \sup_x \| \nat f_{\theta^*} (x)\|(1 - \sigma_-)^{k}+ 4 \sup_x f_{\theta^*} (x) \sum_{i=1}^k (1 - \sigma_-)^{k-i} \|\nat \log q_{\theta^*}(\cdot, \cdot)\|_\infty (1 - \sigma_-)^{i-1}.
\end{align*}
Moreover, using the Fisher identity: 
\begin{align} \label{eq:majoftheta}
&  \sup_x \nat f_{\theta^*}(x)\leq \sup_x f_{\theta^*}(x) \sup_x \nat \log p_{\theta^*}(\Yb_1^n|X_0=x), \nonumber\\
& \quad = \sup_x f_{\theta^*}(x)\sup_x \esp_{\theta^*} \left(\left. \nat \log p_{\theta^*}(\Xb_1^n,\Yb_1^n|X_0=x) \right|\Yb_1^n,\ X_0=x \right),\nonumber \\
& \quad \leq \sup_x f_{\theta^*}(x) \left(\sum_{i=1}^n [\| \nat \log q_{\theta^*}(\cdot,\cdot)\|_\infty+ \| \nat \log g_{\theta^*}(Y_i|\cdot)\|_\infty] \right). 
\end{align}
Combining with (\ref{eq:majoVrais}) completes the proof. 
\end{proof}
As a consequence of Lemma \ref{lem:techniq}, we have
\begin{lemma} \label{lem:techniq2}
  Assume \refhyps{hyp:uniform-ergodicity}{hyp:UniformBound}. Then,  
  \begin{itemize}
  \item [(i)] There exists a random variable $F(\Yb_1^n)$ satisfying $\esps_{\theta^*}(F(\Yb_1^n)^2)< \infty$ and 
$$
\frac{\| \nat \bar p_{\theta^*}(\Yb_1^n)\|}{\inf_x p_{\theta^*}(\Yb_1^n| X_0=x)}\leq F(\Yb_1^n).
$$
\item [(ii)] There exists a constant $C$ such that for all $m,\ k \geq 0$, 
$$
\frac{\|\nat \bar p_{\theta^*}(\Yb_1^n| \Yb_{-k-m}^{-k})-\nat \bar p_{\theta^*}(\Yb_1^n)\|}{\inf_x p_{\theta^*}(\Yb_1^n|X_0=x)} \leq C(E(\Yb_1^n)+ D(\Yb_{-\infty}^{-k})) \lambda^k.
$$ 
  \end{itemize}
\end{lemma}
\begin{proof} As in the proof of Lemma \ref{lem:techniq}, define $f_{\theta}(x):=p_{\theta}(\Yb_1^n|X_0=x)$. Then, 
  \begin{eqnarray*}
    \| \nat \bar p_{\theta^*}(\Yb_1^n)\| &\leq &\left\|\int \nat  f_{\theta^*}(x)  \pi_{\theta^*}(dx)\right\|+ \left\| \int f_{\theta^*}(x)\nat \pi_{\theta^*}(X_{0}=x)\mu(dx)\right\|,\\
&\leq & \sup_{x}  \|\nat  f_{\theta^*}(x)\| + C \sup_{x} f_{\theta^*}(x).
 \end{eqnarray*}
by Lemma \ref{lem:distribInit} (i). Combining with (\ref{eq:majoVrais}) and (\ref{eq:majoftheta}) completes the proof of (i). Now, write: 
   \begin{align*}
&\|\nat \bar p_{\theta^*}(\Yb_1^n| \Yb_{-k-m}^{-k})-\nat \bar p_{\theta^*}(\Yb_1^n)\|\leq \left\|\int \nat p_{\theta^*}(\Yb_1^n|X_{-k})  (\probs_{\theta^*}(dX_{-k}| \Yb_{-k-m}^{-k})-\pi_{\theta^*}(dX_{-k}))\right\|\\
& +\left\| \int p_{\theta^*}(\Yb_1^n|X_{-k}=x)  (\nat\bar p_{\theta^*}(X_{-k}=x| \Yb_{-k-m}^{-k})-\nat \pi_{\theta^*}(X_{-k}=x))\mu(dx)\right\|.\\
  \end{align*}
The first term of the rhs is bounded using Lemma \ref{lem:techniq} (ii), 
$$
\left\|\int \nat p_{\theta^*}(\Yb_1^n|X_{-k})  (\probs_{\theta^*}(dX_{-k}| \Yb_{-k-m}^{-k})-\pi_{\theta^*}(dX_{-k}))\right\| \leq \lambda^k E(\Yb_1^n) \inf_x p_{\theta^*}(\Yb_1^n| X_0=x).
$$
For the second term, fix some $u$ in $\X$. Then, using Lemma \ref{lem:techniq} (i),  
 \begin{align*}
&\left\| \int p_{\theta^*}(\Yb_1^n|X_{-k}=x)  (\nat\bar p_{\theta^*}(X_{-k}=x| \Yb_{-k-m}^{-k})-\nat \pi_{\theta^*}(X_{-k}=x))\mu(dx)\right\|\\
   &\quad =\left\|\int (p_{\theta^*}(\Yb_1^n|X_{-k}=x)-p_{\theta^*}(\Yb_1^n|X_{-k}=u))(\nat\bar p_{\theta^*}(X_{-k}=x| \Yb_{-k-m}^{-k})-\nat \pi_{\theta^*}(X_{-k}=x))\mu(dx)\right\|,\\
& \quad \leq 2 (1-\sigma_-)^{k}(D(\Yb_{-\infty}^{-k})+C)  \left(\frac{\sigma+}{\sigma_-}\right)\inf_x p_{\theta^*}(\Yb_1^n| X_0=x),
 \end{align*} 
which completes the proof of (ii). 
\end{proof}

\noindent \textbf{Proof of Proposition \ref{prop:stationFisher}}. 

\begin{proof} First, write
\begin{align*}
 & \|I_{\Yb_{1}^n}({\theta^*})-I_{\Yb_{1}^n| \Yb_{-k-m}^{-k}}({\theta^*})\|=\sup_{u, \|u\|=1}\left|\esps_{\theta^*} \left\{(u^T \nat \log \denss_{\theta^*}(\Yb_{1}^n))^2-(u^T\nat \log \denss_{\theta^*}(\Yb_{1}^n| \Yb_{-k-m}^{-k}))^2 \right\}\right|,\\
&\quad \leq  \esps_{\theta^*}\left\{\|\nat \log \denss_{\theta^*}(\Yb_{1}^n)-\nat \log \denss_{\theta^*}(\Yb_{1}^n| \Yb_{-k-m}^{-k})\|\|\nat \log \denss_{\theta^*}(\Yb_{1}^n)+\nat \log \denss_{\theta^*}(\Yb_{1}^n| \Yb_{-k-m}^{-k})\|\right\},\\
&\quad \leq\left(\esps_{\theta^*}\left\{\|\nat \log \denss_{\theta^*}(\Yb_{1}^n)-\nat \log \denss_{\theta^*}(\Yb_{1}^n| \Yb_{-k-m}^{-k})\|^2\right\} \right)^{1/2} \\
&\hspace{5cm}\left(\esps_{\theta^*}\left\{\|\nat \log \denss_{\theta^*}(\Yb_{1}^n)+\nat \log \denss_{\theta^*}(\Yb_{1}^n| \Yb_{-k-m}^{-k})\|^2\right\} \right)^{1/2}. 
\end{align*}
It is thus sufficient to prove that 
\begin{align}
  &\lim_{k\to \infty } \sup_{m}\left(\esps_{\theta^*}\|\nat \log \bar p_{\theta^*}(\Yb_1^n| \Yb_{-k-m}^{-k})-\nat \log \bar p_{\theta^*}(\Yb_1^n)\|^2\right)^{1/2}=0,\label{eq:firstone}\\
 &\limsup_{k\to \infty } \sup_{m}\left(\esps_{\theta^*}\|\nat \log \bar p_{\theta^*}(\Yb_1^n| \Yb_{-k-m}^{-k})+\nat \log \bar p_{\theta^*}(\Yb_1^n)\|^2\right)^{1/2}<\infty. \label{eq:secondone} 
\end{align}
We will just prove the first inequality since (\ref{eq:secondone}) is directly implied by (\ref{eq:firstone}). Now,  

\begin{align*}
& \|\nat \log \bar p_{\theta^*}(\Yb_1^n| \Yb_{-k-m}^{-k})-\nat \log \bar p_{\theta^*}(\Yb_1^n)\| \\
&\quad \leq \frac{\|\nat \bar p_{\theta^*}(\Yb_1^n| \Yb_{-k-m}^{-k})-\nat \bar p_{\theta^*}(\Yb_1^n)\|}{\bar p_{\theta^*}(\Yb_1^n| \Yb_{-k-m}^{-k})}+ \frac{\| \nat \bar p_{\theta^*}(\Yb_1^n)\|}{\bar p_{\theta^*}(\Yb_1^n)} \frac{|\bar p_{\theta^*}(\Yb_1^n| \Yb_{-k-m}^{-k})- \bar p_{\theta^*}(\Yb_1^n)|}{\bar p_{\theta^*}(\Yb_1^n| \Yb_{-k-m}^{-k})},\\
&\quad \leq \frac{\|\nat \bar p_{\theta^*}(\Yb_1^n| \Yb_{-k-m}^{-k})-\nat \bar p_{\theta^*}(\Yb_1^n)\|}{\inf_x p_{\theta^*}(\Yb_1^n|X_0=x)}\\
& \quad \quad + \frac{\| \nat \bar p_{\theta^*}(\Yb_1^n)\|}{\bar p_{\theta^*}(\Yb_1^n)} \frac{\sup_{x,x'}|p_{\theta^*}(\Yb_1^n|X_{-k}=x)- p_{\theta^*}(\Yb_1^n|X_{-k}=x')|}{ \inf_x p_{\theta^*}(\Yb_1^n|X_0=x)}, 
\end{align*}
which implies (\ref{eq:firstone}), using Lemma \ref{lem:techniq} (i) and Lemma \ref{lem:techniq2} (i) and (ii). The proof is completed.  
\end{proof}

%\bibliographystyle{ims}
%\bibliography{motherofallbibs}

\end{document}